\documentclass[preprint,12pt]{elsarticle}

\usepackage{amssymb}

\usepackage{UserCommands}

\usepackage{xcolor}
\usepackage{todonotes}
\usepackage{multirow}
\usepackage[official]{eurosym}
\usepackage{soul}
\usepackage{amsmath}
\usepackage[hidelinks]{hyperref}

\makeatletter
\def\ps@pprintTitle{%
	\let\@oddhead\@empty
	\let\@evenhead\@empty
	\def\@oddfoot{\centerline{\thepage}}%
	\let\@evenfoot\@oddfoot}
\makeatother

\begin{document}

\begin{frontmatter}

\title{Convex Optimization of PV-Battery System Sizing and Operation with Non-Linear Loss Models}

\author[inst1,inst2]{Jolien Despeghel}
\author[inst1,inst2]{Jeroen Tant}
\author[inst1,inst2]{Johan Driesen}

\affiliation[inst1]{organization={ESAT-ELECTA},
            addressline={Kasteelpark 10 bus 2445}, 
            city={Leuven},
            postcode={3001}, 
            country={Belgium}}

\affiliation[inst2]{organization={EnergyVille},
            addressline={Thor Park 8310}, 
            city={Genk},
            postcode={3600}, 
            country={Belgium}}

\begin{abstract}
In the literature, when optimizing the sizing and operation of a residential PV system in combination with a battery energy storage system, the efficiency of the battery and the converter is generally assumed constant, which corresponds to a linear loss model that can be readily integrated in an optimization model. However, this assumption does not always represent the impact of the losses accurately. For this reason, an approach is presented that includes non-linear converter and battery loss models by applying convex relaxations to the non-linear constraints. The relaxed convex formulation is equivalent to the original non-linear formulation and can be solved more efficiently. The difference between the optimization model with non-linear loss models and linear loss models is illustrated for a residential DC-coupled PV-battery system. The linear loss model is shown to result in an underestimation of the battery size and cost as well as a lower utilization of the battery. The proposed method is useful to accurately model the impact of losses on the optimal sizing and operation in exchange for a slightly higher computational time compared to linear loss models, though far below that of solving the non-relaxed non-linear problem.
\end{abstract}


\end{frontmatter}

\section{Introduction}

Renewable energy is one of the key solutions to mitigate climate change, by reducing the greenhouse gas emissions of the world's electricity consumption \cite{irena}. On the residential scale, this translates to a rise in distributed generation, with more and more households installing photovoltaic (PV) systems. However, this poses a challenge to the distribution grid. High penetration of PV generation may require grid infrastructure investments in the future \cite{Spilliotis2016}. These investments can be deferred by either limiting the feed-in of distributed generation \cite{Stetz2013} or storing the excess generation by means of PV-battery systems\cite{Moshovel2015}. The latter has the advantage that curtailment of PV generation can be avoided, preventing the loss of renewable generation potential. To maximize this potential, simultaneous optimization of system sizing and operation in combination with more accurate power-dependent modeling of both power converter and battery losses is required. The objective of this paper is to implement more accurate power-dependent loss models in a computationally efficient way in order to assess the need for these loss models in PV-battery optimization and their impact on the optimal sizing and operation, as opposed to linear constant efficiency loss models. 

Past work on the techno-economic sizing and operation of PV-battery systems modeled the power converter and battery losses by means of a constant efficiency, resulting in a linear loss model \cite{Hesse2017,Schram2018,Beck2016,ru2012storage,Linssen2017,Brusco2016,Wu2017, Tant2013}. However, this does not correspond to the real-world creation of losses. For converters, the losses are a non-linear function of the magnitude of the normalized power, and to a lesser extent a function of the voltage \cite{Harrison2008}. For the battery, the losses depend on a multitude of aspects, of which the most important are the (dis)charge power, voltage and state of charge \cite{MousaviG.2014}.  
Few studies have included a more detailed model of the losses in a PV-battery system.

When the losses are modeled in a more detailed manner, increasing the detail of the converter losses is usually prioritized. These losses are represented by means of a non-linear power-dependent function in \cite{Truong2016,Riffonneau2011,Weniger2014,Weniger2016,Ranaweera2016}. In \cite{Truong2016}, 
the battery sizing was specified and the PV size was varied, investigating the profitability of different scenarios.
Where \cite{Riffonneau2011} used a single efficiency curve to model both DC/DC and DC/AC converters, the proposed loss model will differentiate between the DC/DC converter of the battery, the DC/DC converter of the PV and the DC/AC converter.

Besides the converter losses, the battery losses were modeled with more detail using an impedance based model in \cite{Fortenbacher2017,Moshovel2015}. However, due to computational complexity, a simplified model was used instead of the proposed non-linear impedance based model in \cite{Fortenbacher2017}. Note that in \cite{Moshovel2015} the converter efficiency was not only dependent on the power, but also on the voltage.

However, this added complexity is not required to assess the impact of power-dependent loss models on optimal sizing and operation, while negatively impacting the computational burden.

As mentioned above, these non-linear loss models can be implemented to study PV-battery system sizing and operation using one of two computational methods: brute force simulation or mathematical optimization. 
The simulation method was implemented by \cite{Moshovel2015} to study the optimal operation of a PV-battery system by means of a rule-based operational strategy alongside the impact of different management strategies on grid relief. When studying sizing in addition to operation using the simulation method, a range of PV and battery sizes have to be simulated, which was implemented in \cite{Weniger2014} to determine the cost-optimal sizing. Here, the inverter rating was linked in a 1-to-1 ratio to the PV power and the battery capacity. The ratio of the battery inverter to the battery capacity is studied in more detail in later work by Weniger et al. \cite{Weniger2016}, which demonstrates the importance of including battery converter sizing in the optimal sizing of a PV-battery system, given its impact on the battery operation.

An optimization method is generally defined as a method in which a minimization (or maximization) objective and decision variables bound by a set of constraints are defined. When the number of decision variables increases, this method is preferred over simulation as this would lead to a high number of simulations, and thus run time. There are a number of optimization methods that can be used in this context.
Linear optimization and convex optimization were used in \cite{Hesse2017} and \cite{Wu2017}, respectively to cost-optimally size the battery and battery converter. In \cite{Beck2016}, mixed-integer linear programming was used to extend the optimization to include the sizing of the PV system. Furthermore, trade-offs in design decisions can be identified using multi-objective optimization in the context of low voltage distribution grid planning and the potential deferral of grid upgrades \cite{Tant2013}.

The model complexity that arises when implementing more accurate non-linear loss models in optimization can be bypassed by applying model relaxations.
These relaxations, which reduce the model complexity and the resulting run time, have been widely addressed in the field of optimal power flow calculations \cite{Ergun2019}. However, in the field of PV-battery system sizing and operation, the application of relaxations is rather limited. In \cite{Li2016}, second-order cone
programming (SOCP) relaxations were applied to address the non-linearities when sizing distributed energy storage. In \cite{LiConvex2018} the battery model application was limited to the battery charge and discharge scheduling for a limited time horizon. The convexification of non-linear battery models has been illustrated in \cite{LiWen2019}. In \cite{Mehrtash2020}, McCormick relaxations were applied to the non-linearity of battery power flow constraints to size the PV and battery of a zero energy building.

Following the conducted literature review, it has been confirmed that there is limited research available that combines a more accurate quadratic loss model and the optimization of the PV-battery system to the level of the individual converter components, while remaining computationally efficient. However, this topic is important to be able to realize the full potential of PV-battery systems in a cost-efficient way.

The main contributions of this paper are the following: 
\noindent
(1) The development and implementation of accurate quadratic loss models for both the battery and the converters, as opposed to the linear, constant-efficiency loss models generally implemented in literature. Furthermore, the developed loss models are valid across a broad range of battery and converter ratings.

\noindent
(2) The reduction of the high computational burden and complexity associated with non-linear loss models in the original NLP problem by means of convex relaxation, leading to acceptable run times while guaranteeing a single globally optimal solution. Additionally, complementarity constraints are guaranteed without the need for binary variables.

\noindent
(3) The simultaneous optimization of the sizing and operation of the PV-battery system, integrating the impact of the operation on sizing and vice versa by considering the non-linear power-dependent losses. As a result, the optimal sizing of the converters is independent of the PV and battery sizing.

In conclusion, this research has allowed to integrate quadratic loss models using convex relaxation in a computationally tractable way for a time horizon of one year at a resolution of 15 min. By implementing a convex optimization model, a balance is found between the computationally expensive, but exact NLP and the computationally inexpensive LP, which is prone to higher relative errors.

Finally, the paper is structured as follows. Section \ref{section:topology} briefly describes PV-battery system topologies, followed by a more detailed overview of the implemented loss models and their convex relaxations in Section \ref{section:loss_models}. Section \ref{section:optimization} presents the optimization model of the PV-battery system. In section \ref{section:results} the results of the optimization model illustrate the impact of the use of a more detailed loss model. In section \ref{section:discussion} these findings are discussed in a broader context, alongside potential areas of future work. Finally, section \ref{section:conclusion} presents the conclusions. 

\section{PV-battery system topologies}
\label{section:topology}

\begin{figure}
    \centering
    \includegraphics[scale=0.2]{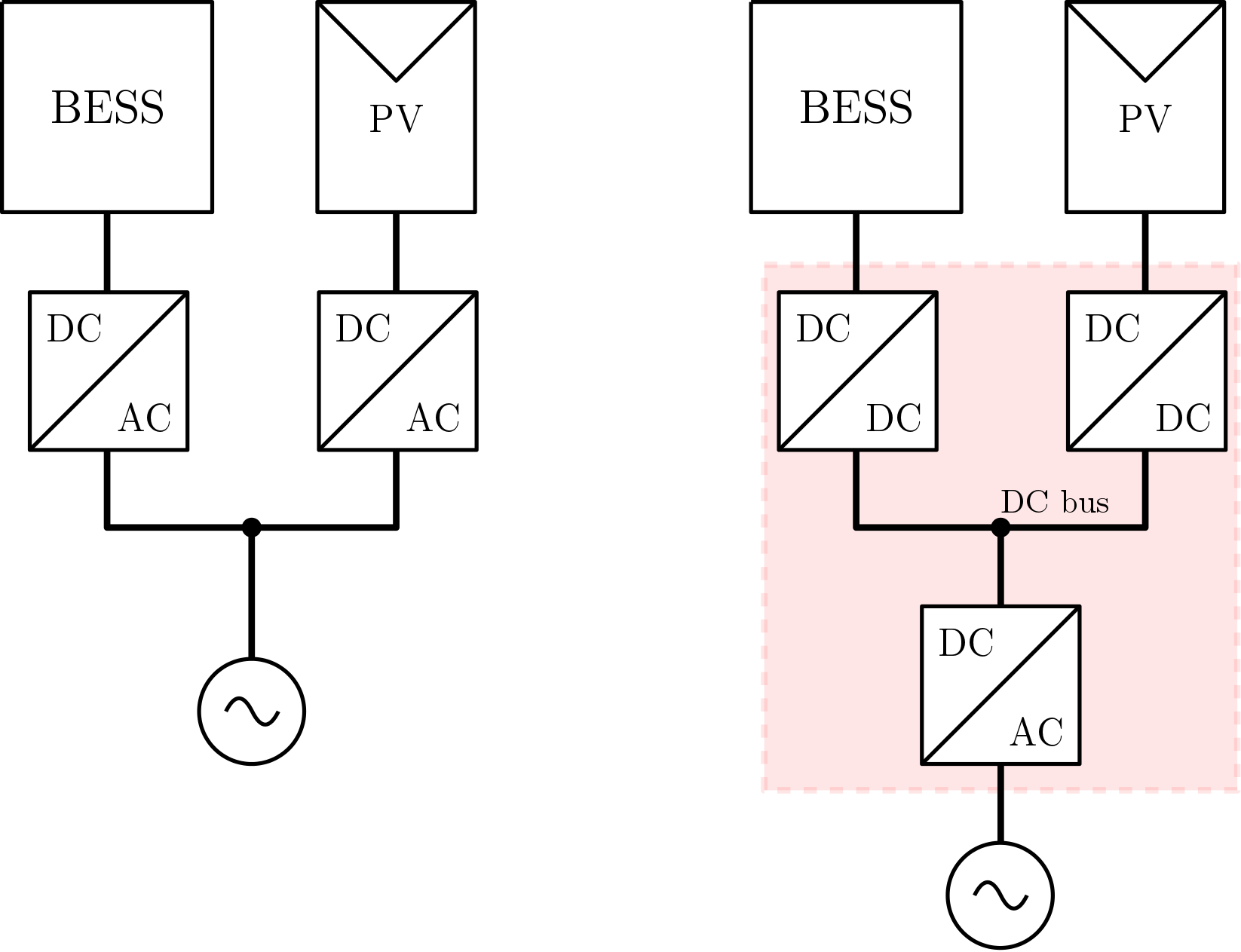}
    \caption{Overview of the two main PV-battery system topologies: AC-coupled (left) and DC-coupled (right)}
    \label{fig:ACvsDC}
\end{figure}

There are two main topologies or configurations of PV-battery systems, as shown in Figure \ref{fig:ACvsDC}: AC-coupled or DC-coupled.
Depending on the system topology, a different converter configuration is implemented. 
The AC-coupled topology requires separate inverters for the PV system and the battery system. The DC-coupled topology requires two separate DC/DC converters, one MPPT tracker and one battery controller, and a single inverter, which are connected to a common DC bus. These converters are often part of the same device, which is then referred to as a hybrid converter.
Regardless of the topology, each of these converters have a certain performance, which can be represented by an efficiency curve, as found on the datasheet, or as a mathematical loss model, fitted on measurement data. 
In the literature, the majority of research has focused on the design and operation of AC-coupled systems \cite{Hesse2017,Schram2018,Beck2016,ru2012storage,Brusco2016,Wu2017,Truong2016,Ranaweera2016,Mehrtash2020,DeOliveiraeSilva2017,VonAppen2015,Weniger2016,Weniger2014}, as opposed to DC-coupled systems \cite{Moshovel2015,Truong2016,Riffonneau2011,Ranaweera2016,Quoilin2016}. The developed models in this work consider both topologies.

\section{PV-battery system loss models and their convex relaxations}
\label{section:loss_models}
This section briefly describes the characteristics of a convex optimization problem before detailing the loss models of the battery and the converters of a PV-battery system, as well as their convex relaxation for use in the optimization model. These loss models are applicable to both DC-coupled and AC-coupled systems.

Following the notation presented in \cite{boyd2004convex}, an optimization problem can be described as follows,
\begin{align}
\min_{x}        &\qquad& f_0(x) & \label{gen_obj}\\
\text{subject to} &      & f_i(x) = 0, & \qquad i= 1,..., m \label{gen_eq_constrait}\\
                 &      & h_i(x) \geq 0, & \qquad i= 1, ...p. \label{gen_ineq_constraint}
\end{align}

 For an optimization problem to be convex both the objective function $f_0(x)$ and the inequality constraint functions $h_i(x)$ need to be convex. The equality constraint functions $f_i(x)$ need to be affine. In these functions, $x$ represents the set of decision variables. The benefit of a convex optimization problem is the guaranteed convergence to the globally optimal solution.
The use of convex relaxations allow to reformulate a non-convex non-linear optimization problem as a convex optimization problem that gives a lower-bound to the solution of the original problem. This technique often leads to the introduction of slack variables to relax non-convex equality constraints as convex inequality constraints. For such relaxations to be exact, the slack variables should be  zero, which can be easily verified by evaluating the original non-linear constraints. 

\subsection{Battery loss model}

Battery loss models in the literature of PV-battery system sizing are often reduced to a single constant efficiency, though the battery efficiency is dependent on many parameters, such as battery state of charge (SOC), battery voltage, charge and discharge current, battery capacity and temperature. 

When deciding on which model to use, a trade-off is made between accuracy and complexity. High accuracy is achieved by electrochemical models, which have the downside of higher complexity due to the use of partial differential equations. Therefore, these models are reserved for battery simulations. A middle ground in terms of both complexity and accuracy is achieved by equivalent circuit models. The lowest complexity, and accompanying low accuracy is found in models where the battery behavior is function of power and energy, and a constant efficiency. To counter the low accuracy of energy-based models, charged-based models have been introduced, though these introduce non-linearities \cite{Aaslid2020}. A range of non-linear battery models and their convex relaxations are presented in \cite{haessig:hal-03032241}. These models include dependencies on energy $E$, power $P$, $P^2$ and $\frac{P^2}{E}$.

In the context of the techno-economic sizing of PV-battery systems, the battery thermal management system can be assumed to work well and no temperature differences are foreseen. In this study one of the main goals is to assess the impact of power-dependent loss modeling, without increasing the modeling complexity to such an extent that it would no longer be tractable. This complexity would occur due to the non-linear nature of the SOC being a fraction of two decision variables. Therefore, this dependency is not included.

The proposed non-linear battery loss model is of the type $\frac{P^2}{E}$, where $E$ corresponds to the rated energy content instead of the instantaneous battery energy. This model allows to determine the loss parameters $\lambda$ and $\gamma$, which are valid for any given battery capacity with the same battery chemistry. First, experimental data were collected by measuring the efficiency of a battery cell at different C-rates $C$ and fitted to  
\begin{align}
    \eta^b = \alpha + \beta C \label{eta_Crate},
\end{align}
 where $ C = \frac{\Pb}{\Ebnom}T^s$ with $T^s$ the sample frequency, the ratio of the battery power $\Pb$ over the nominal energy content $\Ebnom$ on an hourly basis. This ordinary least squares fitting results in the determination of the parameters $\alpha$ and $\beta$.

Second, as the loss is proportional to the complement of the efficiency, the loss model is derived as follows, 
\begin{align}
    \Pbloss &= (1-\eta^b)\Pb \label{loss_init} \\
            &= \underbrace{(1-\alpha)}_{= \lambda}\Pb -\beta \frac{{\Pb}^2}{\Ebnom}T^s & (\alpha > 0) \\
            &= \lambda \Pb -\beta \frac{{\Pb}^2}{\Ebnom}T^s & (\beta < 0) \\
            &= \lambda \Pb +\beta' \frac{{\Pb}^2}{\Ebnom}T^s & (\beta' > 0) \\
    \Pbloss  &= \lambda \Pb + \gamma \frac{{\Pb}^2}{\Ebnom} & (\gamma = \beta'T^s). \label{battery_loss_exact}
\end{align}

Then, \eqref{battery_loss_exact} needs to be relaxed to a convex formulation. The first term is equal to $\Ploss^\text{lin}$ and the second term is equal to $\Ploss^\text{quad}$. 
\begin{align}
    \Ploss^\text{lin} &= \lambda \Pb \label{battery_loss_lin} \\
    \Ploss^\text{quad} &= \gamma \frac{{\Pb}^2}{\Ebnom} \label{battery_loss_quad}
\end{align}
Equality constraint \eqref{battery_loss_quad} is still non-linear and is further reformulated to a rotated second-order cone constraint as follows,
\begin{align}
    \Ploss^\text{quad} \Ebnom \geq \gamma {\Pb}^2. \label{battery_loss_quad_relax}
\end{align}
Finally, in combination with constraints \eqref{battery_loss_lin} and \eqref{battery_loss_quad_relax} the battery loss is completely defined by adding
\begin{align}
    \Ploss &= \Ploss^\text{lin} +\Ploss^\text{quad}.
    \label{battery_loss_relax}
\end{align}

\subsection{Converter loss model}

Converter loss models in the literature of PV-battery system sizing and operation are often reduced to a single constant efficiency \cite{Hesse2017,Schram2018,Beck2016,ru2012storage,Linssen2017,Brusco2016,Wu2017, Tant2013}, which corresponds to a linear loss model. In this loss model the linear coefficient corresponds to the complement of the constant efficiency $\eta$. The loss $\Ploss$ at time step $t$ is then a function of the input power $P_t$:
\begin{align}
    \Ploss = (1-\eta) P_t
    \label{loss_linear}
\end{align}
Though in reality, the losses in a converter are non-linear. 

In the literature these losses have been modeled by means of quadratic models dependent on either the input power \cite{Peippo1994} or the output power \cite{Harrison2008} \cite{El-Aal2006}. A generic quadratic loss model is defined as:
\begin{align}
    \Ploss = \coeffa{} + \coeffb{}P_t + \coeffc{} {P_{t}}^2.  \label{loss_generic}
\end{align}
A similar approach was applied in \cite{Riffonneau2011}, where a quadratic loss model as a function of the normalized input power was fitted. A disadvantage of this model was the use of the efficiency curve of a single inverter to model all converter components, both DC/DC and DC/AC converters, of a DC-coupled system. 
In previous work by the authors \cite{despeghel2019loss}, a quadratic loss model for a DC-coupled PV-battery converter system was proposed and fitted, resulting in separate loss model parameters for the different converter components, namely the DC/AC converter and the DC/DC converters. 

This model was fitted and experimentally validated\cite{despeghel2019loss} for a commercially available DC-coupled converter system, resulting in loss model parameters which are only valid for the rated power of the respective system components. Therefore, the parameters of the loss model need to be transformed in order to preserve the curvature of the quadratic loss model. The converter on which the parameters were fitted can theoretically be placed in parallel $n$ number of times. This leads to the following loss model,
\begin{align}
    \Ploss = n a +b P_t + c \frac{{P_{t}}^2}{n}.
    \label{loss_n}
\end{align}
The loss model parameters $a$, $b$ and $c$ are then transformed by using the original power rating of the converter $P^{nom,og}$ and $n$ is replaced by the desired power rating $\Pnom$ 
\begin{align}
   \Ploss & =  \Pnom \tilde{a} + bP_{t} + \tilde{c}\frac{{P_{t}}^2}{\Pnom}, \label{converter_loss_exact}\\
  &     \text{with } \tilde{a} = \frac{a}{P^{nom,og}} \text{ and } \tilde{c} = c P^{nom,og}  \nonumber.
\end{align}

The integration of the loss model poses a challenge due to the introduction of non-linearities, in this case a quadratic equality constraint. In the field of OPF modeling, quadratic equality constraints are prevalent and require convex relaxations to ensure reasonable run times while guaranteeing optimality \cite{Ergun2019}.

Thus, \eqref{converter_loss_exact} needs to be relaxed to a convex formulation. The first two terms of the equality constraint are grouped in a new linear constraint equal to $\Ploss^\text{lin}$  and the last term is equal to $\Ploss^\text{quad}$. 
\begin{align}
    \Ploss^\text{lin} &= \Pnom \tilde{a} + bP_{t} \label{converter_loss_lin} \\
    \Ploss^\text{quad} &= \tilde{c}\frac{{P_{t}}^2}{\Pnom} \label{converter_loss_quad}
\end{align}
Equality constraint \eqref{converter_loss_quad} is still non-linear and is further reformulated to a rotated second-order cone constraint as follows,
\begin{equation}
    \Ploss^\text{quad} \Pnom \geq \tilde{c} {P_{t}}^2. \label{converter_loss_quad_relax}
\end{equation}
In combination with constraints \eqref{converter_loss_lin} and \eqref{converter_loss_quad_relax} the converter loss is completely defined by adding
\begin{equation}
    \Ploss = \Ploss^\text{lin} +\Ploss^\text{quad}.
    \label{converter_loss_relax}
\end{equation}
The absence of slack in the relaxed inequality constraint presented in \eqref{converter_loss_relax} can be easily determined by evaluating the exact loss equality constraint of \eqref{converter_loss_exact} and comparing the difference. Conversely, the presence of slack would mean energy is shed, increasing the operational cost. This can be countered by means of a cost-minimization objective, as the increase of losses drives up the cost \cite{DeKeyser2017, haessig:hal-03032241}.

This relaxation technique is applied to the loss constraints of the three converters: battery DC/DC converter, PV DC/DC converter and grid inverter.

\section{Optimization Model}
\label{section:optimization}

The following section describes the optimization model of a DC-coupled PV-battery system, of which the loss model constraints were described in Section \ref{section:loss_models}. A similar model can be implemented for AC-coupled PV-battery systems. Additionally, the model can be applied under variable economic assumptions. First an overview of the model is given, followed by the model constraints and considerations on model relaxations and slack.

\subsection{Model Overview}
\begin{figure}
	\centering
	\includegraphics[scale=0.68]{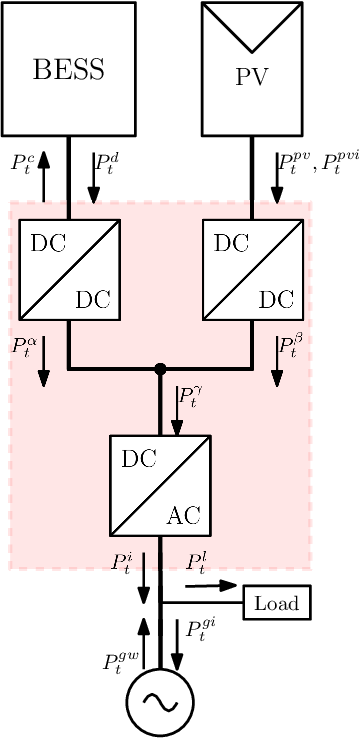}
	\caption{Overview of the DC-coupled PV-battery system and its state variables}
	\label{fig:overview}
\end{figure}

The modeled DC coupled PV-battery system is presented in Figure \ref{fig:overview}. A convex second-order cone programming approach is implemented to optimize the sizing and operation, as the non-linear loss models of the converters and the battery are approximated by convex relaxations. For the purpose of optimal sizing and operation the optimization is performed over all time steps $t$ in a period of one year in 15 minute intervals $\Delta t$.

The main input data consists of the load profile of a single household and the normalized solar generation data of a PV module, which is the DC power at the module terminals expressed in W/Wp. Further technical input data include the parameters of the experimental loss model for the converters of a DC coupled PV-battery system, as detailed in \cite{despeghel2019loss}, and the parameters of the battery loss model. The optimal sizing is found by minimizing the operational and investment costs of the household, as presented in the objective

\begin{multline} 
\mathrm{Min}  \sum_{t} (- \cg{i}\Pgi \dt + \cg{w} \Pgw \dt )T \\
	+ \cdcdc (\Pbnom +\Ppvnom) + \cinv \Pinom  \\ 
	+\cpv \PpvWp + \cbattery \Ebnom .
	\label{primary_objective}
\end{multline} 

with the electricity price $\cg{w}$, feed-in remuneration $\cg{i}$ and investment costs for installed PV $\cpv$, battery $\cbattery$, DC/DC converters $\cdcdc$ and inverter $\cinv$. 
This initial optimization returns the optimal sizing of the PV capacity $\PpvWp$ and the battery capacity $\Ebnom$, as well as the optimal sizing of their respective DC/DC converters $\Ppvnom$ and $\Pbnom$, and of the inverter $\Pinom$. These component sizings are optimal for a certain objective, in this case, the minimization of the consumer's total cost of ownership over a horizon of T years.

Besides the optimal sizing, the model implicitly returns the optimal operation of the system, but for this initial optimization a non-negligible amount of slack remains on the resulting power flows. To reduce this slack and tighten the model constraints a second round of optimization is applied, which will be discussed in more detail in \ref{section:slack}.

The optimal operation, which can function as a benchmark for rule-based control strategies and forecasting methods,  is defined by the following state variables for every time step.

 \begin{itemize}
     \item \Pc, \Pd : the battery charge and discharge power, flowing between the battery and the DC/DC converter
     \item \Ppv : the dc power of the PV modules
     \item \Ppvi: the dc power of the PV generation, possibly curtailed, flowing into the DC/DC converter
     \item \Palph : the power of the battery branch at the DC bus, flowing into or out of the DC/DC converter of the battery
     \item \Pbeta : the power of the PV branch at the DC bus, flowing out of the DC/DC converter
     \item \Pgamma : the net power at the DC bus, flowing into the inverter
     \item \Pinv: the power at the inverter output
     \item \Pgi: the power injected into the grid
     \item \Pgw: the power withdrawn from the grid 
 \end{itemize}

\subsection{PV constraints}

The power output of the PV modules $\Ppv$ is the installed capacity $\PpvWp$ times the normalized solar generation $G^{\text{s}}_t$. In PV systems the available dc power at the PV module can be curtailed, when the generation exceeds the nominal power of the converter. This results in the PV module no longer operating in its maximum power point. In the model this operation mode is included by defining the input power of the DC/DC converter $\Ppvi$, which is allowed to be lower than the available power $\Ppv$  at the module output, while staying below the limit of converter's nominal power $\Ppvnom$.

\begin{align}
	\Ppv &= G^{\text{s}}_t \PpvWp  \label{eq:begin_pf}\\
	\Ppvi &\leq \Ppv \label{PV_curtailment}\\
	\Ppvi &\leq \Ppvnom \label{PV_nominal}
\end{align}

The PV power at the DC bus $\Pbeta$ is defined as follows,
\begin{align}
	\Pbeta = \Ppvi - \Ppvloss, \label{PV_converter}
\end{align}
where $\Ppvloss$ is the loss in the DC/DC converter as defined by \eqref{converter_loss_lin}, \eqref{converter_loss_quad_relax} and \eqref{converter_loss_relax}.

\subsection{Battery constraints}

The battery energy $\Eb$ follows from the energy balance at every time step  $t$,
\begin{align}
	\Eb = \Ebt{t-1} + (\Pc-\Pbcloss) \dt  - (\Pd +\Pbdloss) \dt, \label{Batt_energy}
\end{align}
where the battery losses $\Pbcloss$ and $\Pbdloss$ are defined by \eqref{battery_loss_lin}\eqref{battery_loss_quad_relax}\eqref{battery_loss_relax}.
The battery energy on the first time step is assumed to be equal to the battery energy on the last time step. Furthermore, the charge and discharge power, respectively $\Pc$ and $\Pd$ remain within the bounds of the rated power of the battery converter $\Pbnom$.
\begin{equation}
	0 \leq \Pc \leq \Pbnom   \text{ and }  0 \leq \Pd \leq \Pbnom 
\end{equation}
The battery power at the DC bus $\Palph$ is defined as follows,
\begin{align}
	\Palph = \Pd - \Pdloss - (\Pc - \Pdloss), \label{battery_converter}
\end{align}
where $\Pcloss$ and $\Pdloss$ are the charge and discharge losses of the DC/DC converter.
A more hidden relaxation is the absence of a complementarity constraint for the battery charge and discharge power, which avoids simultaneous charge and discharge of the battery. In \cite{Arroyo2020} counterexamples are shown where this convex storage model relaxation does not hold in the context of MIP-based power system planning and operation problems. For this reason the results are verified to check whether the  battery complementarity condition holds.
\subsection{Grid constraints}

The inverter connects the DC bus to the grid. The inverter output power $\Pinv$ is limited by the nominal power of the inverter $\Pinom$.

\begin{equation}
-\Pinom \leq \Pinv \leq \Pinom
\end{equation}

\subsection{Power balance constraints}
The PV-battery system model is governed by a number of balancing equations, in particular, power balances. These are briefly summarized starting at the PV system, followed by the battery system, DC bus and finally, the AC grid connection.
\begin{align}
\Ppvi - \Ppvloss - \Pbeta &= 0\\
\Pc + \Pcloss -(\Pd - \Pdloss) + \Palph &= 0\\
\Palph + \Pbeta &= \Pgamma \\
\Pgamma - \Pinv - \Piloss &= 0 \\
\Pinv - P^\text{l}_t -\Pgi + \Pgw &= 0 \label{eq:end_pf}
\end{align}

\subsection{Model relaxations and slack} 
\label{section:slack}
Besides the optimal sizing, the model also returns the optimal operation of the system. However, for the objective of  \eqref{primary_objective} a non-negligible amount of slack remains on the resulting power flows. 
Specifically, the battery complementarity is not respected and an overestimation of the losses occurs. To reduce this slack and tighten the model constraints a second round of optimization is applied to minimize the losses. The optimal sizing of the first round of optimization serves as the input and a loss-minimizing objective in combination with a constraint on the objective value is used to further tighten the optimal operation. 
\begin{equation}
\mathrm{Min}  \sum_{t} \Ppvloss+\Pbcloss + \Pbdloss+\Pcloss + \Pdloss + \Piloss 
	\label{objective_finetune}
\end{equation}

The achieved reduction of the power flow slack, results in the optimal operation with sufficiently tight relaxations.

\section{Results}
\label{section:results}

This section presents the results of the optimization for different formulations of the converter and battery loss models, which were obtained on an i7 CPU @2.80 GHz with 32 GB RAM. The model is implemented in Julia \cite{Julia-2017} and uses Gurobi \cite{gurobi} as a solver. Each model formulation will return the optimal sizing and operation while minimizing the total cost of ownership (CAPEX and OPEX). First, the input data of the optimizations are described. Second, the optimal solution of the non-linear loss model formulation is compared to the optimal solution of the convex loss model formulation in order to establish the equivalence of both formulations. Third, the optimal solution of the convex loss model formulation is compared to the optimal solution of the linear loss model formulation as well as combinations of both formulations.

\subsection{Input data}
The optimization model requires a load profile, a PV profile and a set of input parameters. The load profile is obtained from the LINEAR project dataset \cite{Labeeuw2013}, which contains the measured household power consumption with a resolution of 15 min, which corresponds to 33600 time steps. The selected profile has an annual consumption of 2774 kWh. 
The normalized PV generation profile corresponds to a South-oriented PV panel under an inclination of 45 degrees, resulting in an annual available generation of 1020 Wh/Wp.

An overview of the input parameters is presented in Table \ref{tab:param}.
\begin{table}[]
    \centering
    \caption{Overview of parameters used for the optimization of the PV-battery system}
    \begin{tabular}{l c}
    \hline
         Parameter & Value  \\
         \hline
         PV cost $\cpv$ (\euro{}/kWp) & 750\\
         Battery cost $\cbattery$ (\euro{}/kWh) & 250 \\
         DC/DC converter cost $\cdcdc$ (\euro{}/kVA) & 130 \\
         DC/AC converter cost $\cinv$ (\euro{}/kVA) & 200 \\
         grid withdrawal cost $\cg{w}$ (\euro{}/kWh) & 0.26 \\
        grid injection remuneration $\cg{i}$ (\euro{}/kWh) & 0.1 \\
      Investment horizon $T$ (years) & 10\\
\hline         
    \end{tabular}
    
    \label{tab:param}
\end{table}

\subsection{Equivalence of non-linear and convex loss model formulations}
The foundation of the optimization model consists of the power flow constraints, described by \eqref{eq:begin_pf}-\eqref{eq:end_pf}. The objective is either \eqref{primary_objective} for the non-linear model or \eqref{primary_objective} and \eqref{objective_finetune} for the convex model.
The non-linear loss model formulation of the converters in \eqref{converter_loss_exact} and of the battery in \eqref{battery_loss_exact} is compared to the convex loss model formulation, in  \eqref{converter_loss_lin}\eqref{converter_loss_quad_relax}\eqref{converter_loss_relax} and in \eqref{battery_loss_lin}\eqref{battery_loss_quad_relax}\eqref{battery_loss_relax} respectively,  to assess the quality of the convex relaxation as an approximation of the non-linear loss model. The non-linear formulation is the benchmark, representing the physics of reality, but has the disadvantage of a higher numerical complexity. Thus, a convex relaxation is introduced to reduce the computational burden. \\
Due to the high computational complexity of the non-linear formulation the time resolution of the LINEAR load profile as well as the PV profile is increased, reducing the number of time steps. The original profiles represent a year with a resolution of 15 minutes, which corresponds to 33600 time steps. The additional profiles are averaged to 8400 and 16800 time steps, corresponding to a resolution of one hour and 30 min, respectively.
The optimization results of both the original non-linear and relaxed convex formulation are shown in Table \ref{tab:NL_vs_convex}. As both the objective value and the optimal sizing of the convex formulation are identical to that of the non-linear formulation, equivalence between both formulations can be established. The convex formulation is a sufficiently good approximation of the non-linear formulation. Additionally, the run times of the convex formulation with second-order cone programming can be up to 4 times faster compared to the non-linear programming formulation when optimizing for a whole year with a 15 minutes resolution. Note that the run times of non-linear programming generally scale exponentially with increasing problem dimensions, as opposed to the run times of existing polynomial-time algorithms for second-order cone programming.

\begin{table}[!t]

\caption{Optimization Results of Non-linear and Convex Loss Model Formulations}
\label{tab:NL_vs_convex}
\centering

\begin{tabular}{l l l l}
\hline
 & \multirow[t]{2}{*}{time} & non-linear & convex \\
& resolution (h) & & \\
\hline
\multirow{3}{*}{Objective (\euro{})} & 1 &  5107 & 5107 \\
                                & 0.5 & 5292 & 5292 \\
                                & 0.25 & 5494 & 5494 \\
                                \hline
 \multirow{3}{*}{PV size (kWp)}  & 1 & 18.3 &  18.3\\
                                & 0.5 & 17.7 &  17.7\\
                                & 0.25 & 16.7 & 16.7 \\
                                \hline
\multirow{3}{*}{PV DC/DC converter (kVA)} & 1& 11.6 & 11.6  \\
                                & 0.5 & 11.7 &   11.7 \\
                                & 0.25 & 11.6 & 11.6  \\
                                \hline
\multirow{3}{*}{Battery size (kWh)} & 1 & 4.8 & 4.8 \\
                                    &  0.5 & 4.8 &  4.8\\
                    & 0.25 & 4.4 & 4.4 \\
                    \hline
\multirow{3}{*}{Battery DC/DC converter (kVA)} & 1 & 1.2 & 1.2 \\
                            & 0.5 & 1.3 & 1.3 \\
                            & 0.25 & 1.3 & 1.3\\
                            \hline
\multirow{3}{*}{DC/AC converter (kVA)} & 1 & 10.0 & 10.0 \\
            & 0.5 & 10.0 & 10.0 \\
            & 0.25 & 8.5 & 8.5 \\
            \hline
\multirow{3}{*}{Runtime (seconds)} & 1 & 227 & 80  \\ 
                                    & 0.5 & 733 & 393\\
                                    & 0.25 & 2367 & 556 \\
\hline
\end{tabular}
\end{table}

\subsection{Comparison of convex and linear loss model formulations}
The convex loss model is compared to the linear loss model formulation by performing the optimization using the four possible combinations of loss model formulations.
\begin{itemize}
    \item CC-CB: Convex Converter Convex Battery (fully convex)
    \item CC-LB: Convex Converter Linear Battery
    \item LC-CB: Linear Converter Convex Battery
    \item LC-LB: Linear Converter Linear Battery (fully linear)
\end{itemize}

CC-CB serves as the benchmark, as the convex relaxation of the loss model is shown to be a good approximation of the non-linear formulation.
The linear loss model formulations are \eqref{loss_linear} for the converter and \eqref{battery_loss_exact} for the battery. \\

\begin{figure}
	\centering
	\includegraphics[scale=0.6]{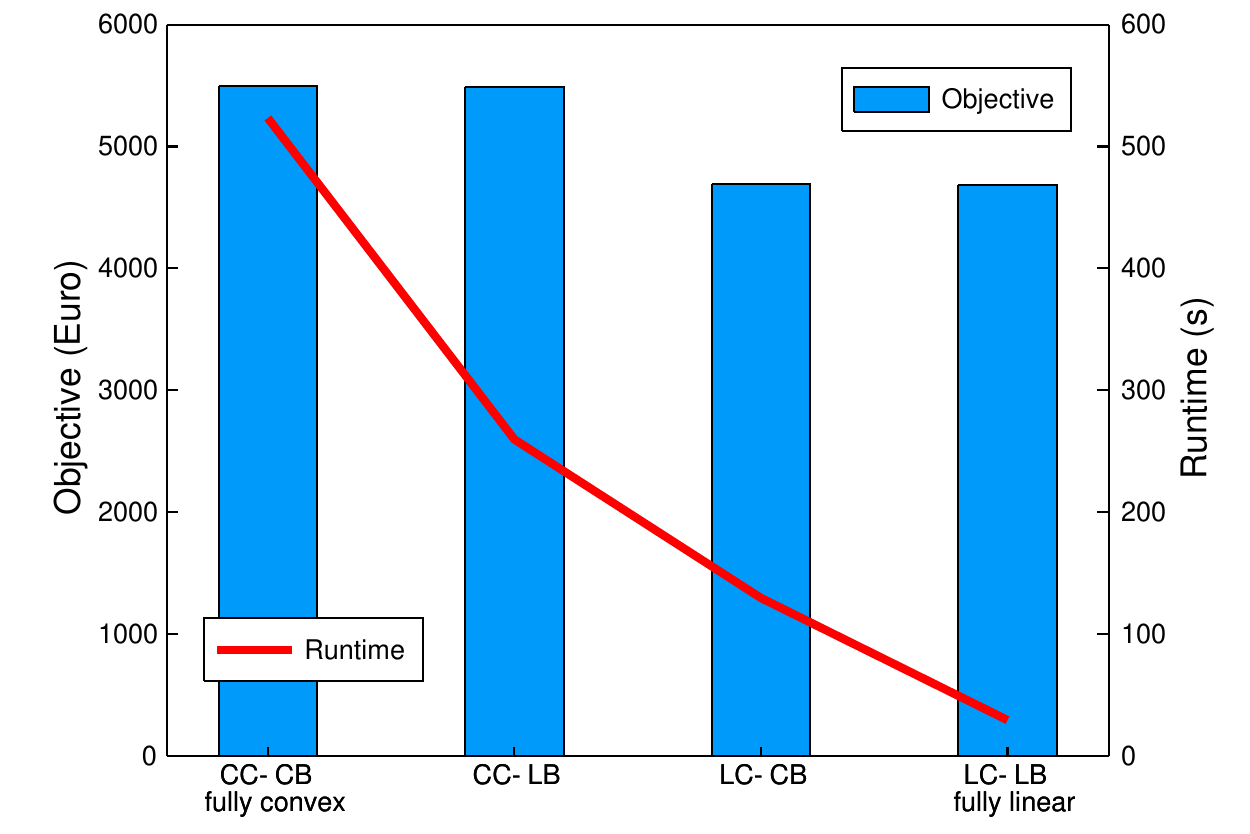}
	\caption{Tradeoff between the runtime and the objective of the different loss model formulations}
	\label{fig:tradeoff}
\end{figure}

The results of the optimization for the four different loss model formulations are shown in Table \ref{tab:LMF_results_overview}. This includes the objective, the system sizing,  as well as the revised objective after convex optimization of the operation for the system sizing of CC-LB, LC-CB and LC-LB. Furthermore, the annual amount of energy injected in and withdrawn from the grid is presented. Lastly, the runtime is presented. The trade-off between the accuracy of the objective and the runtime is illustrated in Figure \ref{fig:tradeoff}.
When the linear loss formulation (LC-LB) is compared to the benchmark (CC-CB), the objective value is underestimated, by 14.7\%. The converter losses have the biggest influence on the underestimation of the objective.

Furthermore, when the operation of CC-LB, LC-CB and LC-LB is revisited using the convex formulation of CC-CB, the objective is no longer underestimated, but slightly overestimated. The overestimation for LC-LB is 1.3\%.

\begin{figure}
	\centering
	\includegraphics[scale=0.6]{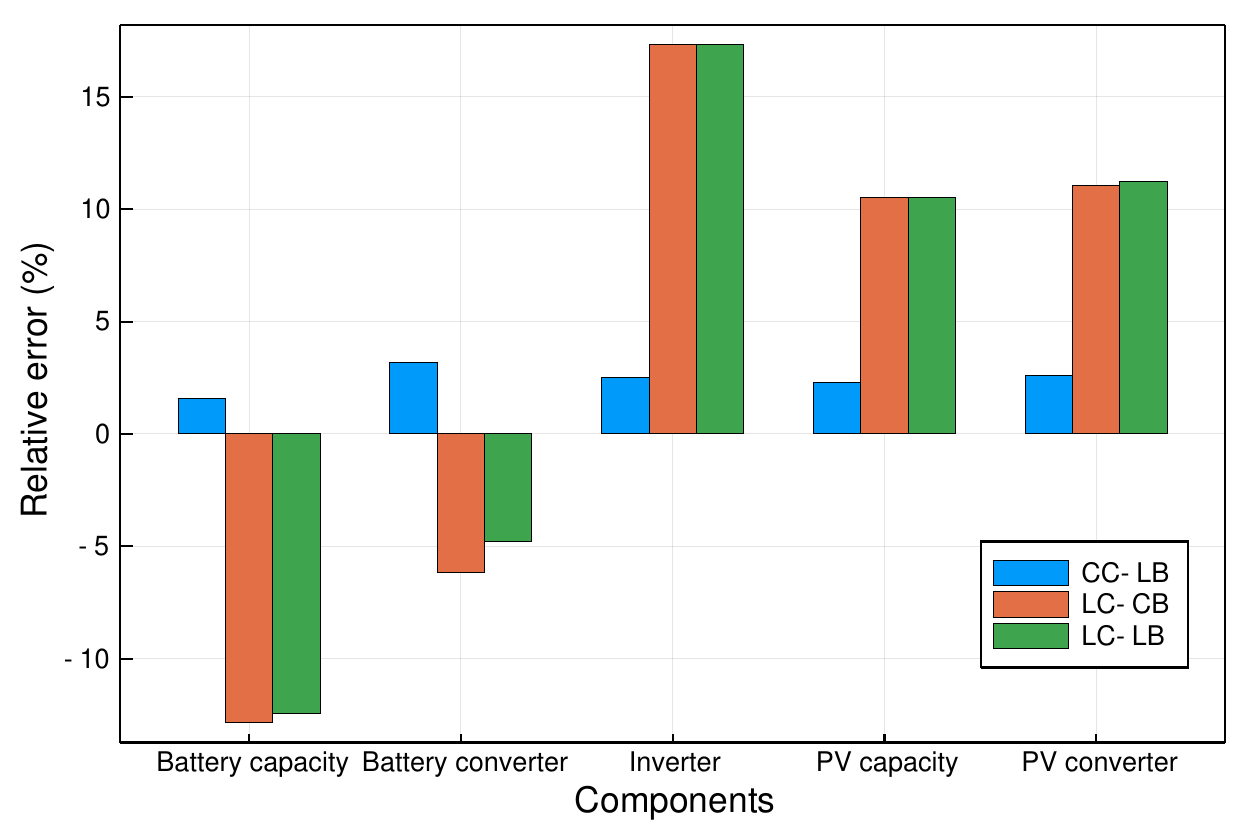}
	\caption{Impact of the loss model formulations on the relative sizing of the PV-battery system}
	\label{fig:LMF_results_sizing}
\end{figure}

Figure \ref{fig:LMF_results_sizing} shows the relative difference in sizing when comparing the benchmark formulation CC-CB to CC-LB, LC-CB and LC-LB. The relative difference in the sizing of the PV-battery system is the largest for the inverter, followed by the battery and the PV components. The inverter rating is overestimated by 17.3\% when optimizing the linear loss model formulation (LC-LB). The battery capacity is underestimated by 12.4\%. Similarly, an underestimation occurs for the battery converter, though smaller, by 4.8\%. The installed PV capacity and the rating of the converter are overestimated by 10.5\% and 11.2\%, respectively. \\

\begin{table}[!t]
\caption{Simulation results for different loss model formulations}
\label{tab:LMF_results_overview}
\centering
\begin{tabular}{l c c c c}
\hline
Loss model & CC-CB  & CC-LB  & LC-CB  & LC-LB \\
            & benchmark &      &       & \\
\hline
 Objective (\euro{}) & 5494 & 5487  &  4693 & 4686\\
 \multirow[t]{2}{*}{Objective under}   & - & 5495 & 5566 &  5565\\
 convex operation (\euro{}) & & & & \\
PV size (kWp)               & 14.8 & 15.2 & 16.4 & 16.4 \\
PV DC/DC converter (kVA)    & 10.2 & 10.5 & 11.3 & 11.4 \\
Battery size (kWh)          & 4.4 &  4.4 & 3.8 & 3.8\\
Battery DC/DC converter (kVA) & 1.3 & 1.4 & 1.3 & 1.3 \\
Grid inverter (kVA)         & 8.5  & 8.7 & 10.0 & 10.0 \\
Grid injection (MWh)        & 11.8 & 12.2 & 13.6 & 13.6\\
Grid withdrawal (kWh)       & 887 & 877 & 732 & 729\\
Runtime (seconds)           & 556 & 268 & 129 & 30\\
\hline
\end{tabular}
\end{table}

\begin{figure}
	\centering
	\includegraphics[scale = 0.6]{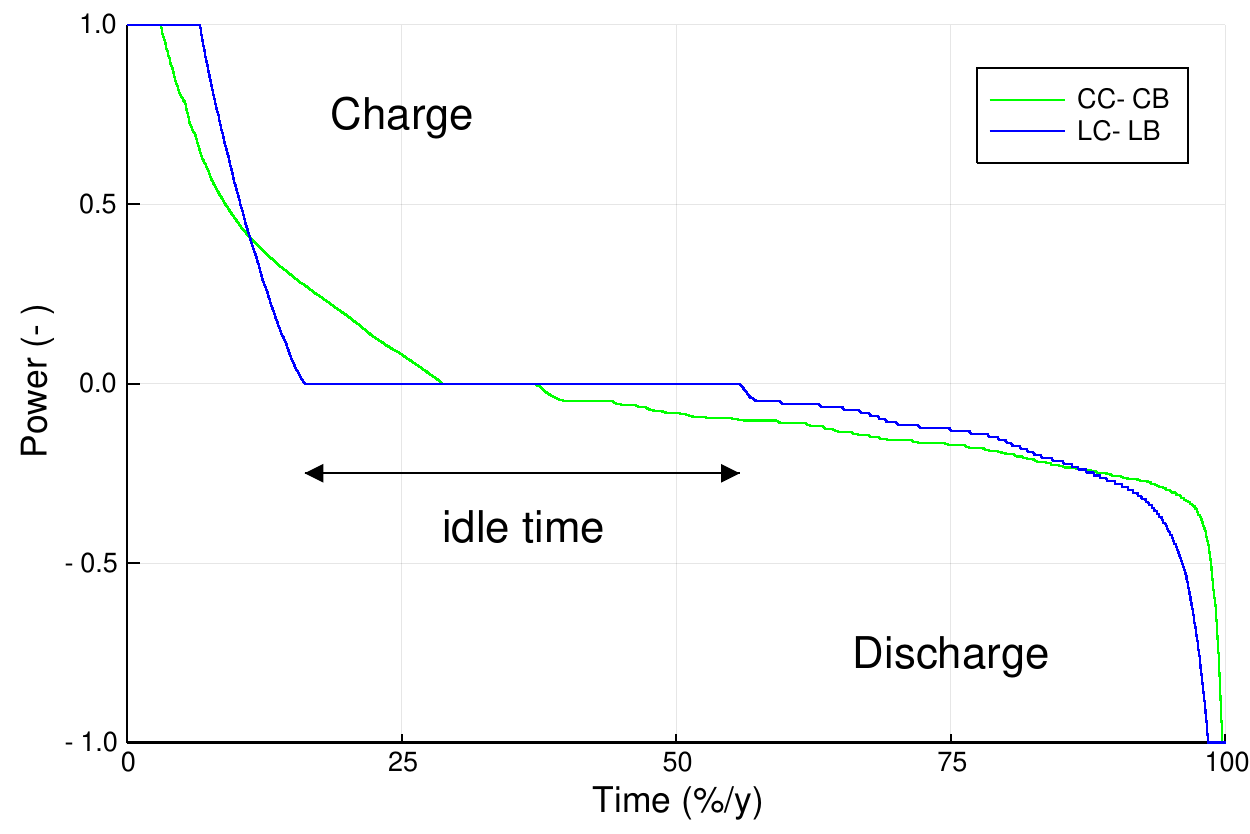}
	\caption{Comparison of the annual battery utilization using normalized battery charge and discharge power duration for CC-CB (blue) and LC-LB (green)}
	\label{fig:results_operation}
\end{figure}

\begin{figure}
    \centering
    \includegraphics[scale = 0.6]{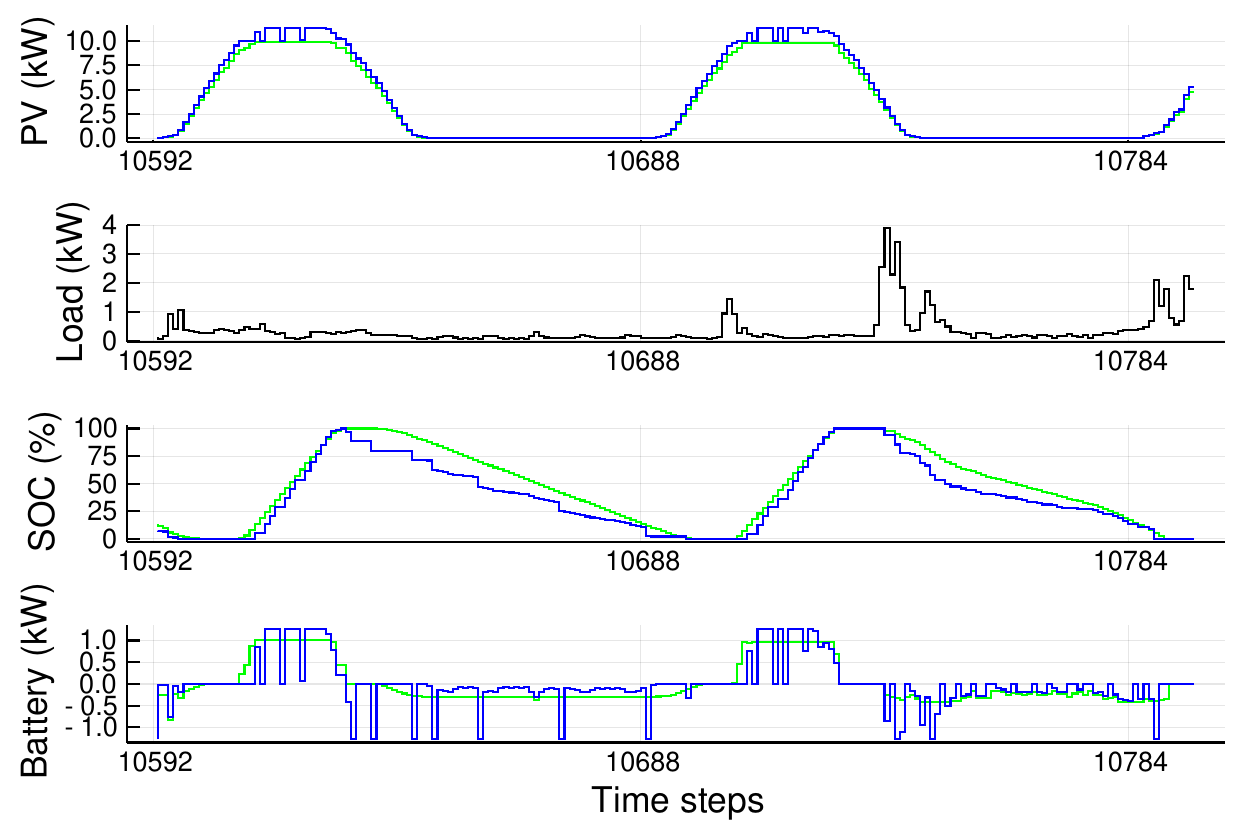}
    \caption{Battery charge (positive) and discharge (negative) power for CC-CB (green) and LC-LB (blue)}
    \label{fig:daily_operation}
\end{figure}

Besides the difference in sizing, the difference in operation between the convex loss model (CC-CB) and the linear loss model (LC-LB) is significant. This is illustrated in Figure \ref{fig:results_operation} for CC-CB and LC-LB,  which displays the normalized battery charge and discharge power duration during a single year. In CC-CB, the battery operates less time at nominal power and more time in the lower power range. In LC-LB, the battery operates more time at nominal power range and is utilized during less hours. This results in an extensive amount of idle time of approximately 40\% compared to 10\%. Additionally, the difference in operation can also be seen on a daily scale. In Figure \ref{fig:daily_operation} the battery charge and discharge powers, as well as the battery state-of-charge (SOC) are shown for two consecutive days. During daytime the battery is charged more smoothly for CC-CB as opposed to LC-LB. After the first day of generation the battery discharges more evenly for CC-CB as compared to LC-LB. On an annual time scale this results in a self-consumption of 20.4\% and and a self-sufficiency of 50.3\% for CC-CB as opposed to 20.6\% and 59.8\%, respectively. The linear loss model thus leads to an overestimation of the self-sufficiency and only a slight overestimation of the self-consumption.

\section{Discussion}
\label{section:discussion}

Different approaches have been used to study the sizing and operation of PV-battery systems. The level of detail that was implemented to model the losses of both the converters and the battery was often limited to a constant efficiency. A constant efficiency is not representative of the non-linear efficiency that has been empirically determined through measurements. Increasing the accuracy of the loss models can assist in realizing the full potential of PV-battery systems, both from a technical and an economical point of view.

In this study, a convex optimization model was presented to optimize the sizing and operation of a PV-battery system. The system losses, of the converters and the battery specifically, were approximated by the convex relaxation of quadratic loss models based on empirical data.
The applied convex relaxations preserved the accuracy of the non-linear optimization result and at the same time reduced the runtime significantly. Consequently, the convex optimization was considered the benchmark for further comparison with the linear optimization. This comparison showed that the linear loss model formulation (LC-LB) can lead to a significant deviation of the model objective and component sizes, though LC-LB could be used to find the initial solution.
Additionally, the convex loss formulation (CC-CB) returns a more balanced operation of the battery, which can have a positive impact on the battery lifetime.\\

Although the findings show that the inclusion of convex loss models allows for a more accurate sizing and operation, these findings are limited to the load profile, PV orientation and economic data of the presented case. Furthermore, the converter loss model was fitted on a single specific inverter, resulting in a loss model that only holds for the specific converter topology, similar to the battery loss model, which was fitted on the measurements of a single type of battery cell. The deviations found between the implementation of convex loss models and linear loss models will be highly dependent on these model parameters.

The focus of future work is twofold. 
On the one hand, including non-linear loss models by means of convex relaxations returns more accurate results. However, the run times are higher. Thus, the conditions need to be assessed for which the use of these convex loss models offer an added value, taking into account the trade-off between level of detail and speed. A further differentiation between the use of a convex converter loss model or battery loss model in combination with a linear loss model can be evaluated, to assess which component is more sensitive to detailed modeling.
On the other hand, the convex optimization model can now be applied to assess a range of policy incentives and economic trends, such as capacity tariffs, grid injection limits, feed-in remunerations and future battery prices. Another future application includes the use of these models in power system operation and planning. Additionally, design considerations such as PV orientation, converter and battery types and system topology can be studied. For the converter, a range of converter topologies could be modeled and selected by means of binary variables. The same could be done for different battery chemistries. Also, the battery degradation could be integrated in the model decision making. Furthermore, the model can be used to evaluate rule-based energy management systems and forecasting strategies.

\section{Conclusion}
\label{section:conclusion}

In this paper, an optimization model has been developed, presented and validated, which includes convex relaxations of non-linear battery and converter loss models. 
These loss models provide a better approximation of the converter and battery behavior compared to linear loss models. At the same time, the runtimes remain far below those of non-linear optimizations, increasing their ease of use.

The presented DC-coupled PV-battery system demonstrated that a more accurate modeling of the battery and converter losses can have a significant impact on the optimization results.
Simultaneously optimizing component sizing and system operation led to a higher utilization and more smooth operation of the battery. Further, for the studied set of model parameters, the linear optimization underestimated the total cost of ownership by 14.7\%.

\section*{CRediT author contribution statement}
\textbf{J. Despeghel: }
Conceptualization, Methodology, Software, Validation, Formal analysis, Writing - Original Draft, Writing - review \& editing, Visualization.
\textbf{J. Tant: }
Conceptualization, Methodology, Writing - Review \& Editing.
\textbf{J. Driesen: }
Resources, Writing - Review \& Editing, Supervision, Funding acquisition

\section*{Declaration of competing interest}
The authors declare that they have no known competing financial interests or personal relationships that could have appeared to influence the work reported in this paper.

\section*{Acknowledgment}
This work is supported by the the KU Leuven Internal Funds for project C24/16/018 and the energy transition funds project BREGILAB organized by the FPS economy, S.M.E.s, Self-employed and Energy.

\bibliographystyle{elsarticle-num-names} 
\bibliography{references.bib}

\end{document}